\documentclass[a4paper,USenglish,numberwithinsect,cleveref, autoref]{lipics-v2021}
\usepackage{mathrsfs}
\usepackage{mathtools}
\hideLIPIcs
\usepackage{mathabx}
\title{F\"uredi--Hajnal and Stanley--Wilf conjectures in higher dimensions}
\author{Yiting Jang}{Universit\'e de Paris, CNRS, IRIF, F-75006, Paris, France \and Department of Mathematics, Zhejiang Normal University, China}{yjiang@irif.fr}{}{}
\author{Jaroslav Ne\v set\v ril}{Computer Science Institute of Charles University (IUUK), Praha, Czech Republic}{nesetril@iuuk.mff.cuni.cz}{https://orcid.org/0000-0002-5133-5586}{}
\author{Patrice Ossona de Mendez}{Centre d'Analyse et de Mathématique Sociales CNRS UMR 8557, France \and Univerzita Karlova v Praze, Czech Republic \and \url{http://cams.ehess.fr/patrice-ossona-de-mendez/} }{pom@ehess.fr}{https://orcid.org/0000-0003-0724-3729}{}
\authorrunning{Y.\ Jiang, J.\ Ne\v{s}e\v{r}il, and P.\ Ossona de Mendez} \newcommand{\ERCagreement}{{\begin{minipage}{.56\textwidth}This paper is part of a project that has received funding from the European Research Council (ERC) under the European Union's Horizon 2020 research and innovation programme (grant agreement No 810115 -- {\sc Dynasnet}). \end{minipage}\hfill\begin{minipage}{.33\textwidth}\includegraphics[width=\textwidth]{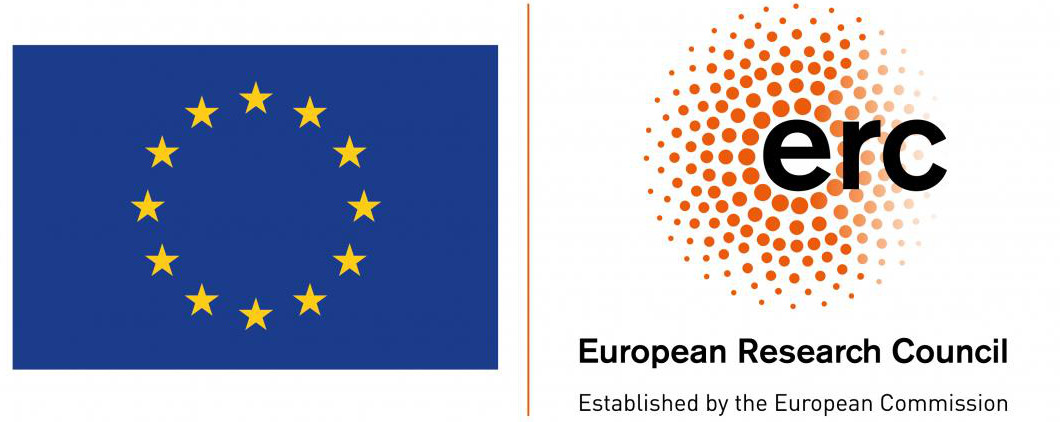}\end{minipage}\hfill}}

\funding{\ERCagreement}
\Copyright{Yiting Jiang, Jaroslav Ne\v{s}e\v{r}il, and Patrice Ossona de Mendez}

\ccsdesc[500]{Mathematics of computing~Combinatorics}
\ccsdesc[500]{Mathematics of computing~Combinatoric problems}
\ccsdesc[300]{Mathematics of computing~Permutations and combinations}

\keywords{binary matrices, grid minors, Latin squares}
\nolinenumbers

\begin{document}

\maketitle
\begin{abstract}
	In this paper we discuss analogs of F\"uredi--Hajnal and Stanley--Wilf conjectures for $t$-dimensional matrices with $t>2$.
\end{abstract}

\section{Introduction}
In~\cite{furedi1992davenport}, F\"uredi and Hajnal considered the problem of determining the maximum of $1$'s in a $0$--$1$ matrix of size $n\times n$ excluding a fixed subconfiguration, and they conjectured that the maximum number $f_2(n,P)$ of  $1$'s in a $0$--$1$ matrix of size $n\times n$ excluding a permutation matrix $P$ as a subconfiguration is bounded by a linear function of $n$.
Klazar \cite{klazar2000furedi} proved that this conjecture, if true, would imply that the Stanley--Wilf conjecture, which states that the growth rate of every proper permutation class is singly exponential. Thus by proving F\"uredi--Hajnal conjecture, Marcus and Tardos~\cite{MarcusT04} also settled Stanley--Wilf conjecture.

In this paper, we are interested in the generalization of both Füredi--Hajnal conjecture and Stanley--Wilf conjecture to higher dimensions, that is in their generalization to the setting of $t$-dimensional matrices (or tensors), with $t>2$. To state these conjectures, we need to introduce several notions.
Let $M$ be an $n_1\times\dots\times n_t$ $t$-dimensional matrix, a typical entry will be denoted as
$M_{i_1,\dots,i_{t}}$.
For $r\in [t]$ and $1\leq j\leq n_r$, the \emph{$(r,j)$-slice} of $M$ is the $n_1\times\dots\times n_{r-1}\times n_{r+1}\times n_t$ $(t-1)$-dimensional matrix $N$ defined by $N_{i_1,\dots,i_{t-1}}=M_{i_1,\dots,i_{r-1},j,i_{r},\dots,i_{t-1}}$. A  \emph{slice} of  $M$ is an $(r,j)$-slice of $M$ for some value of $r$ and $j$.

A \emph{$t$-pattern} $P$ is a $t$-dimensional $0$--$1$ matrix with the property that if $P_{i_1,\dots,i_t}=P_{j_1,\dots,j_t}=1$ then $(i_1,\dots,i_t)$ and $(j_1,\dots,j_t)$ are either equal or differ at at least two positions. A $t$-pattern is \emph{free} if every slice contains at most one entry $1$.
Examples of free $t$-patterns are \emph{$t$-dimensional permutations}, which are the patterns with exactly one entry $1$ in every slice, and examples of non-free $t$-patterns (for $t>2$) are \emph{$t$-dimensional Latin matrices}, which are inductively defined as follows:
$2$-dimensional Latin matrices are the permutation matrices, and, for $t>2$, $t$-dimensional Latin matrices are the patterns with the property that every slice is a $(t-1)$-dimensional Latin matrix. In particular, $t$-dimensional Latin matrices contain $n^{t-1}$ entries $1$, and $3$-dimensional Latin matrices are derived from a standard representation of Latin squares.
	Intermediate between these two extremes are \emph{sunflower} $t$-patterns, which are patterns such that there exist $S\subseteq [t]$ (the \emph{core} of the sunflower $t$-pattern) and $(c_s)_{s\in S}\in [n]^{S}$ with the property that if $P_{i_1,\dots,i_t}=P_{j_1,\dots,j_t}=1$ then for all $s\in S$ we have $i_s=j_s=c_s$ and for all $s\notin S$ we have $i_s\neq j_s$. In particular, $t$-dimensional permutations are sunflower $t$-patterns with empty core.

We say that a $t$-dimensional matrix $M$ \emph{avoids} a $t$-pattern $P$ if $M$ has no submatrix $Q$ with same dimensions as $P$ such that $P_{i_1,\dots,i_t}\leq Q_{i_1,\dots,i_t}$ for all admissible tuple of indices $(i_1,\dots,i_t)$.

For $t\geq 2$, and integer $n$ and a $t$-pattern $P$, 
let $f_t(n,P)$ be the maximum number of entries $1$ in a $t$-dimensional $n\times\dots\times n$ $0$--$1$ matrix avoiding $P$. F\"uredi--Hajnal conjecture stated that for every permutation matrix $P$ (or, equivalently, for every $2$-pattern $P$) we have $f_2(n,P)=O(n)$. We consider here the following generalization of this conjecture.
\begin{conjecture}
The maximum number $f_t(n,P)$ of entries $1$ in a $t$-dimensional $n\times\dots\times n$ $0$--$1$ matrix avoiding a $t$-pattern $P$ is bounded by $a_P\, n^{t-1}$, where $a_P$ depends only on $P$.
\end{conjecture}

Let $T_t(n,P)$ denote the set of all $t$-dimensional $n\times\dots\times n$ $0$--$1$ matrices that avoid the $t$-pattern $P$. Klazar in~\cite{klazar2000furedi} proves that the  F\"uredi--Hajnal conjecture implies Stanley--Wilf conjecture: for every permutation matrix (or, equivalently for every $2$-pattern) $P$, there is a constant $c=c_P$ such that $|T_2(n,P)|\leq c^n$. The natural generalization of this result leads to the following conjecture.

\begin{conjecture}
\label{conj:SWt}
For every $t$-pattern $P$ there exists a constant $c=c_P$ such that the number of $t$-dimensional $n\times\dots\times n$ $0$--$1$ matrices  avoiding  $P$ is such that $|T_2(n,P)|\leq c^{n^{t-1}}$.
\end{conjecture}

The $t$-dimensional Latin matrices defined above were considered in \cite{linial2014upper} as a natural generalization of permutations and Latin squares. In~\cite{linial2014upper}, Linial and Luria prove that the number of $t$-dimensional Latin matrices of order $n$ is bounded by $\bigl((1+o(1))\frac{n}{e^{t-1}}\bigr)^{n^{t-1}}$, and conjecture  that this is actually the correct number. 
Note that Conjecture~\ref{conj:SWt}, if true, immediately implies 
that the number of $t$-dimensional Latin matrices avoiding a $t$-pattern $P$ is bounded by $c^{n^{t-1}}$, where $c$ depends only on $P$, hence the the number of $t$-dimensional Latin matrices would exhibit a jump when a $t$-pattern is excluded, similarly as what happens in the permutation case. 

In this paper we prove these conjectures in  the case where $P$ is a \emph{sunflower} $t$-pattern. As every $2$-pattern is a sunflower $2$-pattern  our results generalizes Marcus and Tardos results. 
Namely, we prove:
\begin{theorem}
\label{thm:FHtg}
For all sunflower $t$-patterns $P$ we have $f_t(n,P)=O(n^{t-1})$.
\end{theorem}
and
\begin{theorem}
\label{thm:SWtg}
The number of $t$-dimensional Latin matrices of order $n$ avoiding a sunflower $t$-pattern $P$ is $O(1)^{n^{t-1}}$.
\end{theorem}

The significance of this last result
appears when considering the (conjectured tight) upper-bound of  
$(n(1+o(1)/e^{t-1})^{n^{t-1}}=O(1)^{n^{t-1}\,\log n}$ for the number of $t$-dimensional Latin matrices~\cite{linial2014upper}.

This paper was also motivated by the recent research related to twin-width of graphs and its possible generalizations to hypergraphs. Due to the space limitations this is only mentioned at the end of this paper.

\section{Full divisions}

As in the proof of Marcus and Tardos, an essential notion is the following.
\begin{definition}
 Let $M$ be an $n_1\times\dots\times n_t$ $t$-dimensional matrix, and let $1\leq k_1\leq n_1$, \dots, $1\leq k_t\leq n_t$.
 A \emph{$k_1\times\dots\times k_t$-division} of $M$ is a tuple
 $(\Pi_1,\dots,\Pi_t)$, where $\Pi_i$ is a partition of $[n_i]$ into $k_i$ intervals. 
 
 The \emph{cell matrix} of  $M$ defined by the $k_1\times\dots\times k_t$-division $\Pi=(\Pi_1,\dots,\Pi_t)$ is the $k_1\times\dots\times k_t$ matrix $M[\Pi]$ with $t$-dimensional entries, such that $M[\Pi]_{i_1,\dots,i_t}$ is the 
 $|I^1_{i_1}|\times\dots\times |I^t_{i_t}|$ $t$-dimensional submatrix of $M$ obtained by keeping only the first indices in the $i_1$th part of $\Pi_1$, \dots, the $t$th index in the $i_t$th part of $\Pi_t$. Entries of $M[\Pi]$ are the \emph{cells} of $M$ in the $k_1\times\dots\times k_t$-division $\Pi$.
 
The \emph{contraction} of $M$ defined by the $k_1\times\dots\times k_t$-division $\Pi$ is the $k_1\times\dots\times k_t$ matrix $N=M/\Pi$ with $N_{i_1,\dots,i_t}=1$ is the cell $M[\Pi]_{i_1,\dots,i_t}$ has at least one non-zero entry, and $0$ otherwise. A $k_1\times\dots\times k_t$-division $\Pi$ of $M$ is \emph{full} if $M/\Pi$ has no zero entries.
\end{definition}
 
 Theorem~\ref{thm:FHtg} will be deduced in Section~\ref{sec:main} from the next theorem.
 \begin{theorem}
 \label{thm:full}
    For every integers $t\geq 2$ and $k$ there exists a constant $\alpha_t(k)$ such that every $n\times\dots\times n$ $t$-dimensional matrix $M$ with more than $\alpha_t(k)\,n^{t-1}$ non-null entries has a full ($t$-dimensional) $k\times\dots\times k$-subdivision.
 \end{theorem}

In this section we proof Theorem~\ref{thm:full} by induction on the dimension $t$, then (for fixed $t$) on induction on $n$. The case $t=2$ has been proved in~\cite{MarcusT04}. Thus we assume $t\geq 3$, and by induction we assume that constants $\alpha_s(k)$ have been defined for every integer $k$ and every $2\leq s<t$. For $2\leq s<t$, we say that an $s$-dimensional $n\times\dots\times n$ matrix is \emph{$k$-heavy} if it contains more than $\alpha_s(k)\,n^{s-1}$ non-null entries.

We define the \emph{$r$th smash} of a $t$-dimensional $n\times\dots\times n$ matrix $M$ as the $(t-1)$-dimensional matrix $P_r(M)$
 such that 
\[
(P_r(M))_{i_1,\dots,i_{t-1}}=\begin{cases}
1&\text{if }\exists 1\leq j\leq n:\ M_{i_1,\dots,i_{r-1},j,i_{r},\dots,i_{t-1}}\neq 0\\
0&\text{otherwise}
\end{cases}
\]

We define an abstract simplical complex $\Delta$, whose maximal faces are the tuples $\{(1,i_1),\dots,(t,i_t)\}$ with $M_{i_1,\dots,i_t}=1$.  It is clear that the number of $t$-faces is the number of entries $1$ in $M$, while the number of $(t-1)$-faces is the sum of the number of entries $1$ in the matrices $P_r(M)$ for $r=1,\dots,t$.

Denote by ${\binom{n}{k}}_t$ the number of $k$-cliques of the Tur\'an graph $T_{n,t}$.
The following lemma is proved in \cite[Theorem 15.1.3]{billera1997face} (see also \cite{frankl1988shadows}).
\begin{lemma}
\label{lem:rep}
Given positive integers $m,k$, and $t$ with $t\geq k$ there are unique $s,n_k,\dots,n_{k-s}$ such that 
\[
m={\binom{n_k}{k}}_t+{\binom{n_{k-1}}{k-1}}_{t-1}+\dots+{\binom{n_{k-s}}{k-s}}_{t-s},
\]
$n_{k-i}-\lfloor\frac{n_{k-i}}{r-i}\rfloor>n_{k-i-1}$ for all $0\leq i<s$, and $n_{k-s}\geq k-s>0$.
\end{lemma}

let ${\rm cl}_i(\Delta)$ denote the number of $i$-faces of a simplicial complex $\Delta$.

\begin{theorem}[\cite{frankl1988shadows}]
\label{thm:shadow}
Let $\Delta$ be an $t$-colorable simplicial complex. Let
\[{\rm cl}_k(\Delta)={\binom{n_k}{k}}_t+{\binom{n_{k-1}}{k-1}}_{t-1}+\dots+{\binom{n_{k-s}}{k-s}}_{t-s}\]
be the unique representation of Lemma~\ref{lem:rep}. Then
\[{\rm cl}_{k+1}(\Delta)\leq {\binom{n_k}{k+1}}_t+{\binom{n_{k-1}}{k}}_{t-1}+\dots+{\binom{n_{k-s}}{k-s+1}}_{t-s}.\]
\end{theorem}

\begin{lemma}
If $t\geq 3$ and $\Delta$ is a $t$-colorable simplicial complex, then
\[{\rm cl}_t(\Delta)\leq 2^t \biggl(\frac{{\rm cl}_{t-1}(\Delta)}{t}\biggr)^{t/(t-1)}.\]
\end{lemma}
\begin{proof}
If ${\rm cl}_{t-1}(\Delta)< t$ then ${\rm cl}_t(\Delta)=0$ and the statement obviously holds. So we assume ${\rm cl}_{t-1}(\Delta)\geq t$.
Note that if $n$ is a multiple of $t$ then $\binom{n}{k}_t=\binom{t}{k}\bigl(\frac{n}{t}\bigr)^k$.
Let 
\[n=t\biggl\lceil \biggl(\frac{{\rm cl}_{t-1}(\Delta)}{t}\biggr)^{1/(t-1)}\biggr\rceil.\]
 Then, according to Theorem~\ref{thm:shadow} we have
\[{\rm cl}_t(\Delta)\leq {\binom{n}{t}}_t
\leq 2^t \biggl(\frac{{\rm cl}_{t-1}(\Delta)}{t}\biggr)^{t/(t-1)}.\] 
\end{proof}

\begin{corollary}
\label{cor:X0}
Let $M$ be a $t$-dimensional matrix and, for $r\in [t]$, let $N_r$ be the number of non-null entries in $P_r(M)$.
Then the number of non-null entries in $M$ is at most  
$2^t\bigl(\frac{1}{t}\sum_r N_r\bigr)^{t/(t-1)}$.
\end{corollary}

\begin{lemma}
\label{lem:col}
Let $r\in[t]$, and 
let $M$ be a $t$-dimensional $\overbrace{p\times\dots\times p}^{r-1}\times n\times\overbrace{p\times\dots\times p}^{n-r}$ matrix
split into $n/p$ $p\times \dots\times p$ blocks $A^1,\dots,A^{n/p}$.
Assume that for more than $(k-1)\binom{p-1}{k-1}^{t-1}$ values $c\in[n/p]$ the $(t-1)$-dimensional matrix $P_r(A^c)$ is $k$-heavy. Then $M$ has a full ($t$-dimensional) $k\times k\times\cdots\times k$-division.
\end{lemma}
\begin{proof}
If $P_r(A^c)$ is $k$-heavy then, by induction hypothesis, $P_r(A^c)$ has a full ($(t-1)$-dimensional) $k\times\dots\times k$-division.
It is known and easy to prove that the number of $k\times\dots\times k$-divisions of a $p\times\dots\times p$ matrix is $\binom{p-1}{k-1}^t$.
Hence, by pigeon-hole principle, for at least $k$ values of $c$, the matrices $P_r(A^c)$ have a common full ($(t-1)$-dimensional) $k\times k\times\cdots\times k$-division. These define a full ($t$-dimensional) $k\times k\times\cdots\times k$-division of $M$.
\end{proof}

We consider a $t$-dimensional  $n\times\dots\times n$-matrix $M$ split into $p\times\dots\times p$ blocks $B^{i,j}$ (with $i,j\in [n/p]$). 
Let $X_r$ ($r\in[t]$) be the set of all  blocks $B^{i,j}$ such that 
the $(t-1)$-dimensional matrix $P_r(B^{i,j})$ is heavy.
According to Lemma~\ref{lem:col}, 
if $A$ has no full ($t$-dimensional) $k\times k\times\cdots\times k$-division
then for every integer $r\in[t]$ we have
\[|X_r|\leq (k-1)\binom{p-1}{k-1}^{t-1}(n/p)^{t-1}.\]
Let $X_0$ be the set of all blocks $B^{i,j}$ that have at least one non-null entry but are not in $X_1\cup\dots\cup X_t$. If $M$ has no full ($t$-dimensional) $k\times\dots\times k$-division then
\[|X_0|\le f_t(n/p,k).\]
Also note that each block $B^{i,j}$ in $X_0$ contains, according to 
Corollary~\ref{cor:X0}, at most $2^t(\alpha_{t-1}(k)p^{t-2})^{t/(t-1)}$ non-null entries.

Therefore 
\[
f_t(n,k)\le \sum_{i=1}^{t}|X_i|p^t+|X_0|2^t(\alpha_{t-1}(k)p^{t-2})^{\frac{t}{t-1}}
\]
Let $c(n,k)=f_t(n,k)/n^{t-1}$. We have
\[
c(n,k)\leq t(k-1)\biggl(\frac{1}{p}\binom{p-1}{k-1}\biggr)^{t-1}+c(n/p,k)2^t\alpha_{t-1}(k)^{\frac{t}{t-1}}p^{-\frac{1}{t-1}}.
\]
Thus if we fix $p=(2\alpha_{t-1}(k))^t$ we get
\[
c(n,k)\leq t(k-1)\biggl(\frac{1}{p}\binom{p-1}{k-1}\biggr)^{t-1}+\frac{c(n/p,k)}{2}.
\]
Thus we can let 
\[
\alpha_t(k)=2t(k-1)\biggl[\frac{1}{2^t\alpha_{t-1}(k)^t}\binom{2^t\alpha_{t-1}(k)^t-1}{k-1}\biggr]^{t-1}.
\]

\section{Proof of the main results}
\label{sec:main}
\begin{theorem}[=Theorem~\ref{thm:FHtg}]
\label{thm:fP}
Let $t\geq 2$ be an integer.
For every sunflower $t$-pattern $P$ we have 
$f_t(n,P)=O(n^{t-1})$.
\end{theorem}
\begin{proof}
As every $2$-pattern is free, the result follows from~\cite{MarcusT04} when $t=2$. So we assume $t\geq 3$.
We prove the statement for all integers $t\geq 3$ by induction on the size of the core of $P$. 

Assume that the core of $P$ is empty, that is that $P$ is a free $k\times\dots\times k$ $t$-pattern.  
Assume for contradiction that there exists an $n\times\dots\times n$ $t$-dimensional $0$--$1$ matrix $M$ avoiding $P$ that has more than $\alpha_t(k)\,n^{t-1}$ entries $1$.
According to Theorem~\ref{thm:full}, the matrix $M$ has a full ($t$-dimensional) $k\times\dots\times k$-division. Selecting an entry $1$ in each cell corresponding to an entry $1$ of $P$ we find a submatrix of $M$ that does not avoid $P$, contradicting the hypothesis. Thus, in this case, $f(n,P)\leq \alpha_t(k)\,n^{t-1}$.

Now assume that the property has been proved (for all integers $t\geq 2$) for all sunflower $t$-patterns with core of size at most $\ell\geq 0$. Let
$P$ be a sunflower $t$-pattern with a core $S$ of size $\ell+1$, and let $(c_s)_{s\in S}$ be the fixed values of the indices corresponding to the positions in the core. Let $s\in S$, let $P'$ be the sunflower $(t-1)$-pattern obtained as the $(s,c_s)$-slice of $P$.
Assume for contradiction that there exists an $n\times\dots\times n$ $t$-dimensional $0$--$1$ matrix $M$ avoiding $P$ that has more than $f_{t_1}(n,P')\,n$ entries $1$. Then there exists $1\leq i\leq n$ such that the $(s,i)$-slice $M'$ of $M$ contains more than $f_{t-1}(n,P')$ entries $1$, which leads to a contradiction as $M'$ avoids $P'$. Thus $f_t(n,P)\leq f_{t-1}(n,P')\,n$. Hence, by induction hypothesis, as the core of $P'$ has size $\ell$, we have $f_t(n,P)\in O(n^{t-1})$.
\end{proof}

Let $T_t(n,P)$ be the number of $t$-dimensional $n\times\dots\times n$ $0$--$1$ matrices avoiding a pattern $P$.
We start by a generalization of the recursion used in Klazar's proof for the $t=2$ case.

\begin{lemma}
\label{lem:TP}
For every pattern $P$ we have
\[
|T_t(2n,P)|\leq |T_t(n,P)|\,(2^{2^t}-1)^{f_t(n,P)}.
\]
\end{lemma}
\begin{proof}
We map $T_t(2n,P)$ to $T_t(n,P)$ by partitioning any matrix $M\in T_t(2n,P)$ into $2\times\dots\times 2$-blocks and replacing each all-zero block by a $0$ entry and all other blocks by $1$-entries. The resulting matrix avoids $P$, and any matrix $N$ in $T_t(n,P)$ is the image of at most $(2^{2^t}-1)^w$ matrices of $T_t(2n,P)$ under the mapping, where $w$ is the number  of $1$-entries of $N$. Here $w\leq f_t(n,P)$ so the recursion follows.
\end{proof}

\begin{theorem}[=Theorem~\ref{thm:SWtg}]
\label{thm:SWt}
The number of $t$-dimensional Latin matrices of order $n$ avoiding a sunflower $t$-pattern $P$ is $O(1)^{n^{t-1}}$.
\end{theorem}
\begin{proof}
According to Theorem~\ref{thm:fP}, for every sunflower $t$-pattern $P$ we have $f_t(n,P)=O(n^{t-1})$. Then it follows from Lemma~\ref{lem:TP} that there exists a constant $c_P$ such that $|T_t(n,P)|\leq c_P^{n^{t-1}}$. As $c_P^{n^{t-1}}$ bounds the number of $n\times\dots\times n$ $t$-dimensional $0$--$1$ matrices avoiding $P$, it bounds (in particular) the number of $t$-dimensional Latin matrices of order $n$ avoiding $P$.
\end{proof}

This result should be compared with the conjecture number of $t$-dimensional Latin squares of order $n$, which is 
$(n(1+o(1)/e^{t-1})^{n^{t-1}}$, thus $O(1)^{n^{t-1}\log n}$.

\section{Remarks}
Marcus--Tardos extremal result is the basic underlying the recently intensively investigated twin-width of graphs (and other binary structures), see~\cite{twin-width1,twin-width2,twin-width3,TWW4_arxiv, TWWP-arxiv, Simon20}.
The results of this paper indicate possibility to generalize the twin-width to hypergraphs (and relational structures of higher arities). Particularly they lead to twin width of $t$-partite $t$-uniform hypergraphs. This can be thought as just the first approximation. In the spirit of Szemer\'edi Regularity Lemma for hypergraphs the defining red edges will have to be considered for all intermediate arities. However for ordered structures this is possible. The details will appear elsewhere.

%\bibliography{latin}

\end{document}